\documentclass[12pt]{amsart}
\usepackage{amsmath,amssymb,amsfonts,amscd,graphicx,xypic}

\usepackage[colorlinks=true, linkcolor=blue, citecolor=blue]{hyperref}

\usepackage{latexsym}
\usepackage{amsmath}

\usepackage{xypic}
\usepackage[all]{xy}

\newtheorem{remark}[subsection]{Remark}
\numberwithin{equation}{section}

\newtheorem{cor}{Corollary}[section]

\newtheorem{lem}[cor]{Lemma}
\newtheorem{theoremsection}[cor]{Theorem}
\newtheorem{conjecturesection}[cor]{Conjecture}
\theoremstyle{definition}

\newcommand{\co}{\colon\thinspace}

\newcommand{\nl}{\hfil\break}

\begin{document}

\title{Tutte relations, TQFT, and planarity of cubic graphs}

\author{Ian Agol and Vyacheslav Krushkal}

\address{Ian Agol\nl
University of California, Berkeley, 970 Evans Hall \#3840, Berkeley, CA, 94720-3840\nl
and\nl
School of Mathematics, Institute for Advanced Study, Einstein Drive, Princeton, NJ 08540 }
\email{ianagol\char 64 berkeley.edu}

\address{Vyacheslav Krushkal\nl Department of Mathematics, University of Virginia,
Charlottesville, VA 22904-4137 USA}
\email{krushkal\char 64 virginia.edu}

\begin{abstract} It has been known since the work of Tutte that the value of the chromatic polynomial of planar triangulations at $(3+\sqrt{5})/2$ has a number of remarkable properties. We investigate to what extent Tutte's relations characterize planar graphs. A version of the Tutte linear relation for the flow polynomial at $(3-\sqrt{5})/2$ is shown to give a planarity criterion for  $3$-connected cubic graphs. A conjecture is formulated that the golden identity for the flow polynomial characterizes planarity of cubic graphs as well. In addition, Tutte's upper bound on the  chromatic polynomial of planar triangulations at $(3+\sqrt{5})/2$ is generalized to other Beraha numbers, and an exponential lower bound is given for the value at  $(3-\sqrt{5})/2$. The proofs of these results rely on the structure of the Temperley-Lieb algebra and more generally on methods of topological quantum field theory.
\end{abstract}

\maketitle

\section{Introduction}

\begin{figure}[ht]
\includegraphics[height=2.2cm]{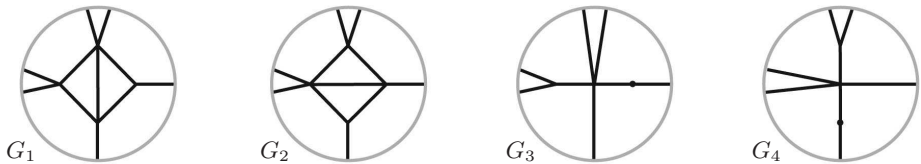}
{\scriptsize
    \put(-348,3){$G_1$}
    \put(-252,3){$G_2$}
    \put(-159,3){$G_3$}
    \put(-66,3){$G_4$}}
\caption{The graphs in the identity (\ref{linear Tutte}).}
\label{fig:linear Tutte}
\end{figure}

In the late 1960s W. T. Tutte observed \cite{T1, T2} that the value of the chromatic polynomial of planar triangulations 
at $(3+\sqrt{5})/2$ obeyes a number of surprising relations.
The Tutte golden identity states
\begin{equation} \label{golden identity}
 {\chi}^{}_T({\phi}+2)=({\phi}+2)\; {\phi}^{3\,V-10}\, ({\chi}^{}_T({\phi}+1))^2,
\end{equation}
where $T$ is any planar triangulation, $V$ is the number of its vertices,  and ${\phi}=(1+\sqrt{5})/2$ denotes the golden ratio.
Tutte also established a linear identity for the chromatic polynomial at ${\phi}+1$:
\begin{equation} \label{linear Tutte}
{\chi}^{}_{G_1}({\phi}+1)+{\chi}^{}_{G_2}({\phi}+1)=
{\phi}^{-3}({\chi}^{}_{G_3}(\phi+1)+{\chi}^{}_{G_4}({\phi}+1)),
\end{equation}
where $G_i$ are planar graphs which are identical outside a disk and are related within the disk 
as shown in figure \ref{fig:linear Tutte}. 

Using duality  between the chromatic polynomial ${\chi}^{}_G$ of a planar graph $G$ and the flow polynomial $F^{}_{G^*}$ of its dual, (\ref{golden identity}) and 
(\ref{linear Tutte}) may be restated as identities for the flow polynomial of planar cubic graphs.

This paper is motivated in part by the question of whether such identities for the chromatic and flow polynomials characterize planarity of graphs. We prove a converse to Tutte's theorem (\ref{linear Tutte}) for $3$-connected cubic graphs: such a graph is planar if and only if a version of the linear relation (\ref{linear Tutte}) for the flow polynomial holds at each edge, 
see Theorem \ref{planar theorem}.
We also consider the identity (\ref{golden identity}) for the chromatic and flow polynomials. 
We conjecture that the golden identity for the flow polynomial characterizes planarity of cubic graphs, see Conjecture \ref{golden conjecture} for a precise statement. On the other hand, examples are given of non-planar graphs satisfying the chromatic polynomial version of the golden identity, and we explain the difference between the chromatic and flow settings in section \ref{inequality section}.

In \cite{FK} the identities (\ref{golden identity}), (\ref{linear Tutte}) were shown to fit in the framework of (2+1) dimensional topological quantum field theory (TQFT).
In particular, these identities were proved to be a 
consequence of the structure of the Temperley-Lieb algebra at roots of unity, and the linear relation (\ref{linear Tutte}) was shown to have a counterpart at each value $2+2\, cos(\frac{2{\pi}j}{n})$.
We extend this framework and define a category for the flow polynomial which takes into account abstract (not necessarily planar) graphs. The structure of this category, together with combinatorial properties of non-planar cubic graphs, give the characterization 
of planarity in Theorem \ref{planar theorem}.

We also use the TQFT framework to generalize Tutte's bound on the value of the chromatic polynomial of a planar triangulations $T$ with $V$ vertices
\begin{equation} \label{inequality}
|{\chi}_T({\phi}+1)|\leq {\phi}^{5-V}
\end{equation}
to all Beraha numbers $B_n=2+2\, cos(\frac{2{\pi}}{n+1})$ and to real values $x\geq 4$, see theorem \ref{Cauchy Schwartz theorem}. The upper bound (\ref{inequality}) is interesting in connection to the conjecture of Beraha that large planar triangulations have a real zero close to ${\phi}+1$. This conjecture remains open. In contrast to (\ref{inequality}), in theorem \ref{lower bound theorem} we give a lower bound at the Galois conjugate value $(3-\sqrt{5})/2$, which is exponentially increasing in $V$.

\section{The Tutte polynomial, the Temperley-Lieb algebra,  and the Jones-Wenzl projectors.} \label{background section}

To provide background for the results of this paper and to fix the notaion, this section recalls basic information about graph polynomials and the Temperley-Lieb algebra, and it summarizes the results of \cite{FK}. 
The Tutte polynomial $T_G(x,y)$ of a graph $G$ satisfies
the contraction-deletion rule: for any edge $e$ of $G$ which is not a bridge or a loop,  $T_G=T_{G\smallsetminus e}+T_{G/e}$. 
The Tutte polynomial of a graph consisting of $b$ bridges and $l$ loops is defined to be $x^b y^l$. This paper will
mainly concern two specializations, the chromatic polynomial ${\chi}_G$ and the flow polynomial $F_G$:
$${\chi}_G(x)= (-1)^{V-c(G)} x^{c(G)} T_G(1-x, 0), \; \; F_G(x)= (-1)^{E+V+c(G)} T_G(0,1-x),$$

where $V, E$ denote the number of vertices and edges of $G$, and $c(G)$ is the number of connected components.
For a planar graph $G$ and its dual $G^*$, one has $T_G(x,y)=T_{G^*}(y,x)$. If $G$ is connected, this implies 
\begin{equation} \label{duality}
F_G(x)=x^{-1} {\chi}_{G^*}(x).
\end{equation}

To describe the TQFT context for Tutte's relations (\ref{golden identity}), (\ref{linear Tutte}), 
we summarize the relevant facts about the Temperley-Lieb algebra, the chromatic algebra and their structure at roots of unity.
The reader is referred to \cite{FK} for more details.

\subsection{The Temperley-Lieb algebra} \label{background subsection}

\begin{figure}[ht]
\includegraphics[height=2cm]{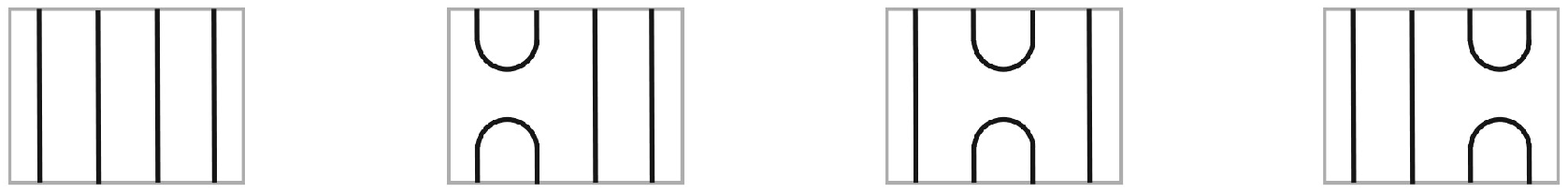}
{\small     \put(-412,23){$1 =$}
    \put(-317,25){$e_1 =\frac{1}{d}$}
     \put(-208,25){$e_2 = \frac{1}{d}$}
    \put(-100,25){$e_3 = \frac{1}{d}$}
    }  
\vspace{.2cm} \caption{Generators of $TL_4$}
\label{fig:TL}
\end{figure}

The Temperley-Lieb algebra in degree $n$, $TL_n$, is an algebra over
${\mathbb C}[d]$ consisting of linear combinations of
$1$-dimensional submanifolds in a rectangle. 
Each submanifold meets
both the top and the bottom of the rectangle in $n$ points. 
The submanifolds are considered equivalent if they are
isotopic rel boundary. Deleting a simple closed
curve has the effect of multiplying the element by $d$. 
When $d$ is specialized to a complex number, the algebra will be denoted
$TL^d_n$.
The multiplication is given by vertical stacking of rectangles. The usual generators
of $TL_4$ are illustrated in figure \ref{fig:TL}.

The trace $tr_d\co TL^d_n\longrightarrow {\mathbb C}$ is defined on
rectangular pictures (additive generators) by connecting the top
and bottom endpoints by disjoint arcs in the complement of the rectangle in the
plane,
and then evaluating $d^{\# circles}$. The Hermitian product 
is defined by $\langle a,b\rangle=tr(a\, \bar b)$, where the
involution $^-$ reflects pictures in a horizontal line and
conjugates the
coefficients.

The {\em trace radical} of $TL_n^d$ is the ideal consisting of
the elements $a$ such that $tr(ab)=0$ for all $b\in TL^d_n$.
An important fact, underlying the construction of SU(2) quantum invariants of knots and $3$-manifolds, is that 
the trace radical in the Temperley-Lieb algebra is
non-trivial precisely at the special values $d\, = \, 2\, \cos \left(\frac{\pi j}{n+1}\right)$. (In fact, these are the values of $d$ for which there exists a  non-trivial ideal in the
Temperley-Lieb category, and the ideal equals the trace radical.)
Moreover, for these values of $d$, (assuming $j$, $n+1$ are coprime)
the trace radical is generated by a special element, the {\em Jones-Wenzl
projector}  $P^{(n)}$.  At primitive roots of unity ($j=1$) the Hermitian product descends to a 
positive definite inner product on the quotient of $TL^d_n$ by the trace radical.

\subsection{The chromatic algebra} 
\begin{figure}[ht]
\includegraphics[bb=20 0 539.256 76.32, height=1.85cm]{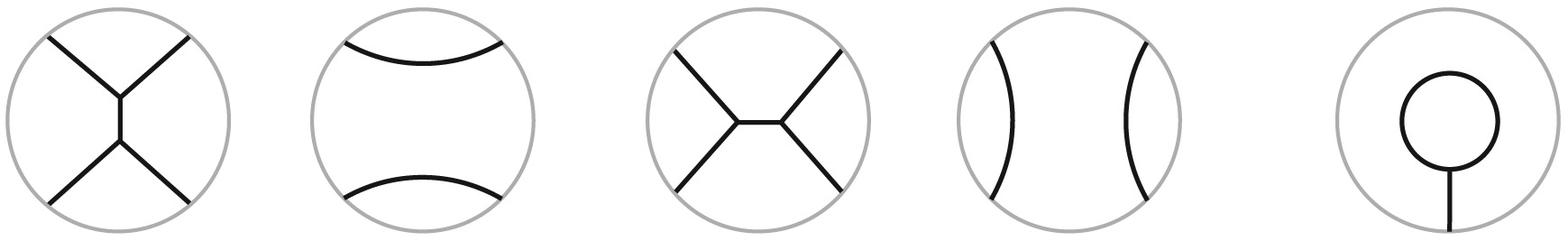}
    \put(-318,23){$+$}
    \put(-247,23){$=$}
    \put(-175,23){$+$}
    \put(-23,23){$=0.$}
    \put(-111,8){$,$}
\caption{Defining relations in the the chromatic algebra.}
\label{chromatic relations}
\end{figure}

${\mathcal C}_n$ is an algebra over ${\mathbb C}[Q]$, 
whose elements are formal linear combinations of 
isotopy classes of planar cubic graphs in a
rectangle $R$, modulo the local relations shown in figure \ref{chromatic relations}.
(The first relation is the version  to the contraction-deletion rule for cubic graphs.)

\begin{figure}[ht]
\includegraphics[height=1.85cm]{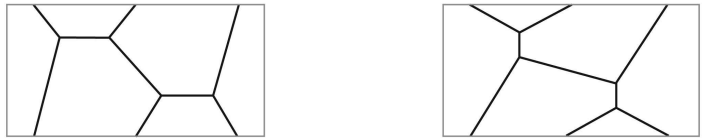}
\caption{Examples of graphs in the definition of ${\mathcal C}_3$.}
\label{ChromaticGen}
\end{figure}

The intersection of a graph with the boundary
of the rectangle $R$ consists of $2n$ points: $n$ points at the top and
the bottom each, figure \ref{ChromaticGen}. 
It is convenient to allow $2$-valent vertices as well, 
and the value of a
simple closed curve is set to be $Q-1$. When $Q$ is set to be a specific complex number, 
the algebra is denoted ${\mathcal C}^Q_n$.

The trace $tr\co {\mathcal C}_n^Q\longrightarrow{\mathbb C}$ is defined by connecting the
endpoints of $G$ by disjoint arcs in the complement of the rectangle $R$ in the plane
and evaluating the flow polynomial of the resulting graph at $Q$. (Equivalently, the trace equals $Q^{-1}$ times the chromatic 
polynomial of the dual graph. The focus in \cite{FK} was on the chromatic polynomial, explaining the name of the algebra.)
For example, the traces of the elements in figure \ref{ChromaticGen} are $0$, $(Q-1)(Q-2)^2$. 
The trace is well-defined since the local relations are precisely the relations defining the flow polynomial of a planar cubic graph.
The multiplication and the Hermitian product on ${\mathcal C}^Q_n$ are defined analogously to the case of the Temperley-Lieb algebra.

\subsection{The map ${\mathbf{{\mathcal C}_n\longrightarrow TL_{2n}}}$} \label{algebra map}
\begin{figure}[h]
\includegraphics[height=1.8cm]{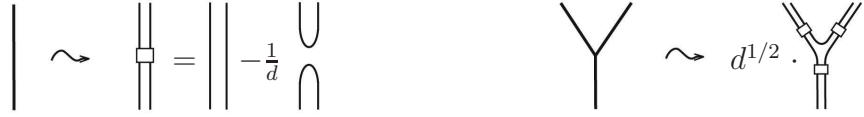}
    \put(-265,21){$=$}
    \put(-240,21){$-\frac{1}{d}$}
    \put(-54,21){$d^{1/2}\;\cdot$}
\caption{The homomorphism ${\Phi}\co {\mathcal C}^Q_n\longrightarrow TL^d_{2n}$, $Q=d^2$.}
\label{fig:map}
\end{figure}

Consider the homomorphism ${\Phi}\co {\mathcal C}^Q_n\longrightarrow
TL^d_{2n}$, where $Q=d^2$. This map replaces 
each  edge of a graph  with the second Jones-Wenzl projector,  
and 
each vertex is resolved as shown in figure \ref{fig:map}. Moreover, for a trivalent graph $G$ there is an overall factor
$d^{V/2}$, where $V$ is the number of vertices of $G$.

One checks that ${\Phi}$ induces a well-defined homomorphism of algebras ${\mathcal C}^Q_n\longrightarrow TL^d_{2n}$,
where $Q=d^2$, and moreover it is trace-preserving: the diagram
\begin{equation} \label{chromaticTL}
\xymatrix{ {\mathcal C}_n^{{d^2}}  \ar[d]^{tr} \ar[r]^{\Phi} & TL_{2n}^{d}  \ar[d]^{tr_{d}}\\
{\mathbb C}\ar[r]^=  &  {\mathbb C} }\end{equation}
commutes.
It follows that the pullback under ${\Phi}$
of the trace radical in $TL^d$ is in the trace radical of ${\mathcal
C}^{d^2}$. Note that elements in the trace radical in the chromatic algebra 
are precisely local relations on graphs which preserve the flow polynomial, 
or equivalently the chromatic polynomial of the dual graph.
When $d={\phi}$, the trace radical of the Temperley-Lieb algebra is
generated by the Jones-Wenzl projector $P^{(4)}$.   The Tutte relation (\ref{linear Tutte}) then
may be seen as a consequence of the structure of $TL^{\phi}$ since it is mapped by $\Phi$ to
$P^{(4)}$ \cite{FK}.

\section{A planarity criterion from TQFT} \label{planarity section}

\begin{figure}[ht]
\includegraphics[height=2.1cm]{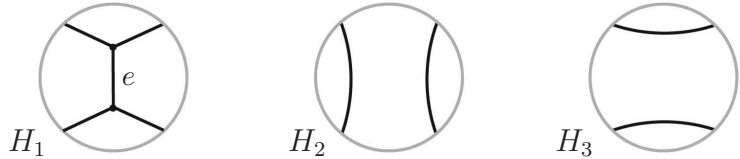}
{\small \put(-236,27){$e$}}
    \put(-279,2){$H_1$}
    \put(-173,2){$H_2$}
    \put(-72,2){$H_3$}
\caption{The graphs $H_i$ in the identities (\ref{3 term}), (\ref{ConjugateLinearRelation}). The graphs are identical outside the disk.}
\label{three graphs}
\end{figure}

It follows from duality (\ref{duality}) that Tutte's linear relation (\ref{linear Tutte}) has an immediate analogue, a 4-term identity for the flow polynomial of dual graphs. Further, using the relations in figure \ref{chromatic relations} (equivalent to the contraction-deletion rule), it can be restated as a 3-term identity 
\begin{equation} \label{3 term}
 {\phi}{F}^{}_{H_1}({\phi}+1)\, - \,  {F}^{}_{H_2}({\phi}+1)\, +\, {\phi}^{-1}\, {F}^{}_{H_3}({\phi}+1)\, =\, 0,
 \end{equation} 

 where planar graphs $H_i$ are shown in figure \ref{three graphs}.

We will use a version of this identity at the Galois conjugate value $(3-\sqrt 5 )/2$:

\begin{lem} \label{cubic Tutte} {\sl 
Consider three planar graphs $H_i$ which are locally related as shown in figure \ref{three graphs}. Then the values of their flow polynomials at 
${\phi}^{-2}=(3-\sqrt 5 )/2$ satisfy the linear identity}
\begin{equation} \label{ConjugateLinearRelation}
 {F}^{}_{H_1}({\phi}^{-2})\, +\, {\phi} \,  {F}^{}_{H_2}({\phi}^{-2})\, +\, {\phi}^2\, {F}^{}_{H_3}({\phi}^{-2})\; =\; 0.
\end{equation}
\end{lem}

Recall from section \ref{background section} that the identity (\ref{3 term}) holds since its left hand side is the pull-back of the $4$th 
Jones-Wenzl projector $P^{(4)}$ from the Temperley-Lieb algebra $TL^{d}$ to the chromatic algebra ${\mathcal C}^Q$. Here $d={\phi}$, $P^{(4)}$ generates the trace radical, and the corresponding value of $Q$ is  $d^2={\phi}+1$. The identity (\ref{ConjugateLinearRelation}) holds for the same reason, being a pull-back of $P^{(4)}$ for $d=(1-\sqrt{5})/2=-{\phi}^{-1}$ and $Q=(3-\sqrt{5})/2={\phi}^{-2}$.

Consider the question of whether (\ref{ConjugateLinearRelation}) holds for a {\em ribbon} (not necessarily planar) graph. Recall that a ribbon graph is an abstract graph $G$ together with an embedding into a surface $S$ so that the complement $S\smallsetminus G$ consists of $2$-cells. 
Given such an embedding, a small neighborhood $N$ of $G$ in $S$ is a compact surface with boundary, which may be thought of as a choice of a $2$-dimensional thickening of $G$. (Such a thickening of vertices and edges may be alternatively encoded using cyclic ordering of half-edges incident to every vertex.)  For our applications it does not matter whether the surface is orientable, and in fact only a local ribbon structure will be used, where a thickening of a given edge and of adjacent half-edges is specified.

\begin{figure}[ht]
\includegraphics[height=2cm]{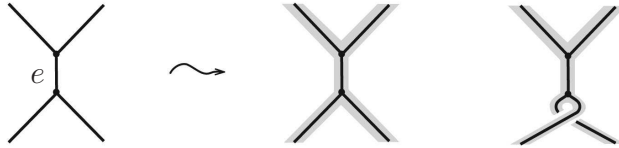}
\put(-227,25){$e$}
\caption{Two ribbon structures near an edge of a cubic graph.}
\label{fig:ribbon}
\end{figure}

A ribbon structure is needed to make sense of the relation (\ref{ConjugateLinearRelation}).
Specifically, an edge $e$ of an abstract cubic graph $H_1$, shown on the left in figure \ref{three graphs}, uniquely determines $H_3$. However the graph denoted $H_2$ in the figure depends on a choice of a ribbon structure near the edge $e$.
There are two such choices, figure \ref{fig:ribbon}.

We are in a position to state the main result of this section:

\begin{theoremsection} \label{planar theorem} \sl A $3$-connected cubic graph $G$ is planar if and only if  the identity (\ref{ConjugateLinearRelation}) for the flow polynomial at $(3-\sqrt{5})/2$ holds for at least one of the two ribbon structures at every edge of $G$.
\end{theoremsection}

The theorem follows from the following lemma. 
Note that in the first part of the statement the graph is not assumed to be cubic.

\begin{lem} \label{key lemma} \sl \nl
(1) Let $G$ be a graph with two adjacent trivalent vertices. Let $e$ denote the edge connecting them, and fix one of the two ribbon structures near $e$ (figure \ref{fig:ribbon}). Denote by $V, E$ the number of vertices, respectively edges, of $G$. Then  
\begin{equation} \label{lemma inequality}
(-1)^{V-E}[ {F}^{}_{H_1}({\phi}^{-2})\, +\, {\phi} \,  {F}^{}_{H_2}({\phi}^{-2})\, +\, {\phi}^2\, {F}^{}_{H_3}({\phi}^{-2}) ] \geq 0,
\end{equation}

where $H_1=G$ and $H_2,H_3$ are locally related to $H_1$ as in figure \ref{three graphs}. 

(2) 
Moreover, if $G$ is a $3$-connected cubic graph then $G$ is non-planar if and only if there exists an edge $e$ with a strict inequality in (\ref{lemma inequality}) for both ribbon structures at $e$. 
\end{lem}

The proofs of lemma \ref{key lemma} and theorem \ref{planar theorem} will be given in the context of {\em flow category} which is introduced next.

\subsection{The flow category} 
Abstract (not necessarily planar) graphs may be used to define an algebra along the lines of the Temperley-Lieb and chromatic algebras
in section \ref{background section}. However since 
the algebra structure will not be needed in the proof of theorem  \ref{planar theorem}, and also the graphs are no longer assumed to be planar, 
the notion of a category is more suitable.

\begin{figure}[ht]
\begin{center}
\includegraphics[height=1.9cm]{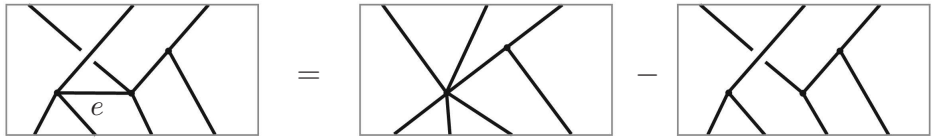}
\put(-242,22){$=$}
\put(-114,22){$-$}
{\small \put(-320,10){$e$}}
\caption{The contraction-deletion rule in ${\mathcal F}_{4,3}$.}
\label{FlowGen}
\end{center}
\end{figure}

The objects of the flow category are finite ordered sets $\overline{n}=\{ 1,\ldots, n\}$. 
Consider ${\mathcal G}^{}_{m,n}=\{$finite graphs with $m+n$  marked univalent vertices$\}$, where the marked vertices are divided into two ordered subsets of $m$, respectively $n$ vertices. The edges incident to the marked vertices
are called {\em boundary} edges and the rest are {\em internal} edges. 
Morphisms between $\overline{m}$, $\overline{n}$ in the flow category are elements of ${\mathcal F}^Q_{m,n}$: formal ${\mathbb C}$-linear
combinations of graphs ${\mathcal G}^{}_{m,n}$, modulo the contraction-deletion relation (figure \ref{FlowGen}) which applies to internal edges.
In addition, the value of a loop is set to be $Q-1$, and graphs which have a univalent vertex (other than the specified marked vertices) 
are set to be zero. It follows from the defining relations that a $2$-valent vertex incident to an interior edge can be removed from a graph.

Graphs whose equivalence classes are elements of ${\mathcal F}^Q_{m,n}$ may be represented geometrically as in figure \ref{FlowGen}.
(Over/under-crossings do not carry any information here since the figure is meant to represent only the abstract graph structure and not a specific planar projection.)

A graph without marked vertices (and so no boundary edges), considered in ${\mathcal F}^Q_{0,0}\cong {\mathbb C}$, evaluates to its flow polynomial at $Q$.
The pairing ${\mathcal F}^Q_{k,m}\times {\mathcal F}^Q_{m,n}\longrightarrow {\mathcal F}^Q_{k,n}$ is obtained by gluing along $m$  boundary edges.
In particular, this pairing applied to two graphs $A\in {\mathcal G}^{}_{0,m}$, $B\in {\mathcal G}^{}_{m,0}$ gives $\langle A, B\rangle=$ the value of the flow polynomial 
$F^{}_{A\cup B}$ at $Q$.

\begin{remark} \rm 
The flow category may be extended to a definition of a TQFT. This is a special case of  a well-known construction \cite{Wa} of TQFTs as fields (``pictures'') on a manifold modulo local relations. In particular, the example of the Tutte polynomial and of the local relation given by the contraction-deletion rule was discussed in \cite{Wa1}. Unlike the Temperley-Lieb and the chromatic algebra which give rise to 
the (2+1)-dimensional (SU(2) and SO(3) respectively) theories, this TQFT is not specific to dimension $2$ since it is formulated in terms of abstract, not embedded graphs. 
We use elements of the structure of the flow category as a convenient setting for the proof of theorem \ref{planar theorem}.
\end{remark}

\begin{figure}[ht]
\begin{center}
\includegraphics[height=2cm]{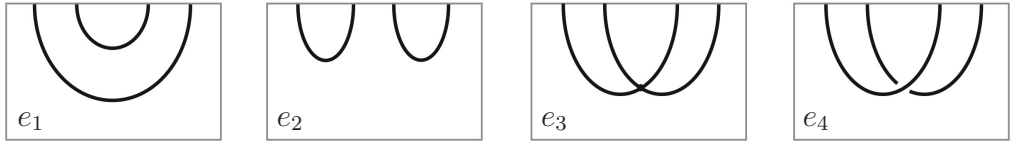}
    \put(-378,7){$e_1$}
    \put(-280,7){$e_2$}
    \put(-180,7){$e_3$}
    \put(-81,7){$e_4$}
\caption{A basis of ${\mathcal F}^Q_{0,4}$.}
\label{fig:basis}
\end{center}
\end{figure}

Given any graph representing an element of ${\mathcal F}^Q_{m,n}$,  the contraction-deletion rule may be used to eliminate all internal edges. In particular, four graphs $e_i$ in figure \ref{fig:basis} form a basis of  ${\mathcal F}^Q_{0,4}$. Note that three of these graphs (viewed relative to a fixed embedding of the marked vertices in the boundary of a rectangle) are planar and one, denoted $e_4$, is non-planar.

{\em Proof of lemma  \ref{key lemma}.}
Let $G$ be a graph with two adjacent trivalent vertices (connected by an edge $e$),  as in the statement of the lemma, and fix $Q={\phi}^{-2}=(3-\sqrt{5})/2$. 
Cut each of the four edges adjacent to $e$ in half, and consider the graph $G$ as the union $\overline H_1\cup \overline G$. Here  $\overline H_1$ is a ``neighborhood'' of the edge $e$, that is the part of $G$ shown on the left in figure \ref{three graphs},  and $\overline G$ is the rest of the graph $G$. $\overline H_1$ is considered as an element of ${\mathcal F}^Q_{0,4}$ and $\overline G$ an element of ${\mathcal F}^Q_{4,0}$, and the union $G=\overline H_1\cup \overline G$ is taken along $4$ pairs of marked points. 

Consider the element $P\in {\mathcal F}^Q_{0,4}$,
$$
P:=  \overline H_1 \, +\, {\phi} \,  \overline H_2 \, +\, {\phi}^2\, \overline H_3.
$$
where $\overline H_i$ denotes the part of the graph $H_i$ contained in the disk in figure \ref{three graphs}, $i=1,2,3$.
(In the planar context, this is an element of the trace radical of the chromatic algebra ${\mathcal C}^Q$,  section \ref{algebra map}.)
Lemma \ref{key lemma} may then be restated as saying that 
\begin{equation} \label{restatement}
(-1)^{V-E}\langle P, \overline G\rangle \geq 0.
\end{equation}

The proof of (\ref{restatement}) is by induction on the number of internal edges of $\overline G$. 
First consider the base of the induction: no internal edges. Since the first three basis vectors in figure \ref{fig:basis} are planar and 
$P$ is in the trace radical of the chromatic algebra,  $\langle P, e_i \rangle =0$, $i=1,2,3.$  A direct calculation gives $\langle P, e_4 \rangle =-{\phi}$, so $(-1)^{V-E}\langle P, e_4 \rangle >0$.

For the inductive step, consider the contraction-deletion rule
\begin{equation} \label{inductive}  (-1)^{V-E}\langle P, G' \rangle =  (-1)^{V-E}  \langle P, G'/  e' \rangle - (-1)^{V-E}  \langle P,  G'\smallsetminus  e' \rangle,
\end{equation}

where $ e'$ is an edge (neither a  bridge nor a loop) of $G'$, and $V, E$ denote the number of vertices and edges of $H_1\cup  G'$. 
Both $  G'/  e'$, $ G'\smallsetminus  e'$ have fewer edges than $ G'$, so the inductive hypothesis applies. 
Moreover, the sign of $(-1)^{V-E}$ of $  G'/ e'$ agrees with that of $ G'$, and it differs from the sign for
$ G'\smallsetminus  e'$. Therefore $(-1)^{V-E}\langle P, G' \rangle \geq 0$, being the sum of two non-negative terms.

To conclude the inductive step, we need to analyze the effect of removing  a loop from $G'$. This changes the sign $(-1)^{V-E}$, but the loop value is $-{\phi}^{-1}$, so overall removing a loop
does not affect the sign of $(-1)^{V-E}\langle P, \overline G \rangle$. This completes the proof of part (1) of the lemma.

\begin{figure}[ht]
\begin{center}
\includegraphics[height=2.4cm]{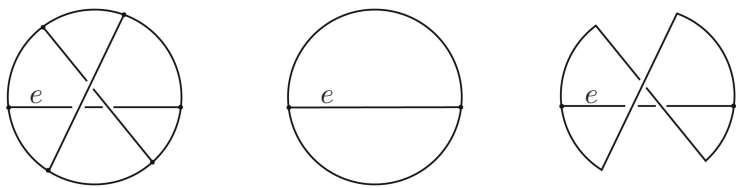}
{\small    
\put(-268,32){$e$}  
\put(-158,32){$e$} 
\put(-58,32){$e$} 
}
\caption{Left: a subdivision of
$K_{3,3}$ with three pairwise non-adjacent non-subdivided edges (drawn straight). The two figures on the right show 
the two relevant minors.}
\label{fig:k33}
\end{center}
\end{figure}

One of the implications in the second statement in lemma \ref{key lemma} is clear: if $G$ is planar, then  the equality in (\ref{lemma inequality}) follows from 
lemma \ref{cubic Tutte}.
The proof of the other implication relies on the theorem \cite{Ke, Th} that any $3$-connected non-planar graph with at least six
vertices contains a cycle with three pairwise crossing chords. In other words, the theorem
asserts that any $3$-connected non-planar graph, different from $K_5$, contains a subdivision of
$K_{3,3}$ with three pairwise non-adjacent non-subdivided edges, figure \ref{fig:k33}.

We will use only a weaker version of the theorem, that one of the edges is unsubdivided, call it $e$.
As in the proof above of part (1) of the lemma, decompose $G=\overline H_1 \cup \overline G$. 
The two ribbon structures near $e$ correspond to the two cyclic orderings of the four univalent vertices of
$\overline H_1$ (and correspondingly of $\overline G$). 
As illustrated in figure \ref{fig:k33}, for either cyclic ordering $\overline G$ has the non-planar basis element $e_4$ (figure \ref{fig:basis}) as a minor.

Fix either ribbon structure at $e$.
Using the contraction-deletion rule on the internal edges of $\overline G$, $(-1)^{V-E}\langle P, \overline G\rangle$ equals the sum of 
$2^k$ terms, where $k$ is the number of internal edges. Each term is of the form $(-1)^{V-E}\langle P, e_i\rangle$ 
where $e_i$ is a basis element in figure \ref{fig:basis}. Moreover, it is shown in the proof of (1) above that each of the $2^k$
summands is non-negative. As indicated in the previous paragraph, at least one of the terms involves the non-trivial pairing with $e_4$.
It follows that $(-1)^{V-E}\langle P, \overline G\rangle>0$.
\qed

\section{Bounds for the chromatic and flow polynomials and the Beraha numbers} \label{bounds section}

It has been known since the work of Tutte \cite{T1, T2} and Beraha \cite{Beraha} that  the Beraha numbers $B_n=2+2\, cos(\frac{2{\pi}}{n+1})$, where $n$ is a positive integer, play a special role in the theory of the chromatic polynomial of planar graphs. 

\begin{theoremsection} \label{Cauchy Schwartz theorem} \sl
Given a planar triangulation $T$,  let $x$ be either a Beraha numbers $B_n$ or a real number $\geq 4$. Then  
\begin{equation} \label{Cauchy Schwartz} 
|{\chi}_T(x)|\; \leq \; x(x-1)(x-2)^{(V-2)}.
\end{equation}
\end{theoremsection}

For $x=B_4= {\phi}+1=(3+\sqrt{5})/2$ this upper bound coincides with the Tutte estimate
\begin{equation} \label{Tutte inequality}
|{\chi}_T({\phi}+1)|\leq {\phi}^{5-V}.
\end{equation}

 Note that the upper bound (\ref{Cauchy Schwartz}) is exponentially decreasing only when $x-2<1$. 
 It is an interesting question whether the bound may be improved for 
$4$-connected planar triangulations.

{\em Proof of theorem \ref{Cauchy Schwartz theorem}.}
When a given triangulation $T$ has a multiedge
or a triangle not bounding a face, one may decompose into
two triangulations and induct to get the inequality. 
More precisely, the chromatic polynomial of a graph $G$ which is the union of two subgraphs $G_1, G_2$ such that 
$G_1\cap G_2$ is the complete graph $K_n$, satisfies
$${\chi}^{}_G(x)\, =\, \frac{{\chi}^{}_{G_1}(x)\, {\chi}^{}_{G_2}(x)}{{\chi}^{}_{K_n}(x)}.$$

For planar graphs the relevant cases are $n=2,3$, and if both $G_1, G_2$ satisfy (\ref{Cauchy Schwartz}), then 
so does $G$.
Therefore it suffices to prove the statement for
$4$-connected planar triangulations $T$.
By a theorem of Whitney \cite{Whitney} $T$ has
a Hamiltonian cycle. The cycle divides $T$ into
two outerplanar triangulations. The doubles of these triangulations
have chromatic polynomial the same as that of the outerplanar
graph, which is  $x(x-1)(x-2)^{(V-2)}$. As explained in section \ref{algebra map}, the chromatic polynomial of a planar graph may be
evaluated in the Temperley-Lieb algebra. The inner product in $TL^d$ is positive semi-definite 
for $d\in \{2\, cos(\frac{\pi}{n+1})\}$, and it is positive definite for $d>2$ \cite{Jo}, \cite[Corollary 7.2]{FK2}.
Then (\ref{Cauchy Schwartz}) follows from the Cauchy-Schwarz inequality.
\qed

\begin{theoremsection} \label{lower bound theorem} \sl
Given a planar triangulation $T$, its chromatic polynomial satisfies 
\begin{equation} \label{conjugate bound}
|{\chi}^{}_T({\phi}^{-2})| \, \geq \, {\phi}^{2V-6}.
 \end{equation}
\end{theoremsection}

{\em Proof.} 
Consider the $3$-connected cubic planar graph $G=T^*$, then  $|{\chi}^{}_T({\phi}^{-2})|=\, {\phi}^{-2}\, |F^{}_G ({\phi}^{-2}) |$.
It follows from an Euler characteristic calculation that the numbers $n, V$ of vertices of $G, T$ are related by $n=2V-4$. Then  (\ref{conjugate bound}) is equivalent
to the inequality
\begin{equation} \label{lower bound}
 |F^{}_G ({\phi}^{-2}) | \; \geq \; {\phi}^{n-2}.
\end{equation}

for the flow polynomial of $3$-connected planar cubic graphs $G$ with $n$ vertices.
Consider the operation of {\em removing an edge} $e$ from $G$, where one disregards the $2$-valent vertices (the boundary of the removed edge) in $G\smallsetminus e$, so that the resulting graph is again cubic.  In this terminology the graph $H_3$ in figure \ref{three graphs} is obtained by removing $e$ from $H_1$.
A theorem of Barnette - Gr${\rm \ddot{u}}$nbaum \cite{BG} (and a closely related result of Tutte \cite{T66}) states that for any $3$-connected graph with more than $6$ edges there exists an edge whose removal yields a graph which is also $3$-connected. Said differently, there is a sequence of edge removals so that each graph in the sequence is $3$-connected and the last graph  is $K_4$.

Consider an inductive application of the identity (\ref{ConjugateLinearRelation}) to such a sequence of edge removals. Let the graphs in figure \ref{three graphs} show the inductive step: the edge $e$ is removed from a planar $3$-connected cubic graph $H_1$, giving rise to a planar $3$-connected cubic graph $H_3$. In particular, it follows that all three graphs $H_i$ are connected.

The flow polynomial of a connected bridgeless graph with $n$ vertices and $E$ edges is non-zero, of sign $(-1)^{(E-n+1)}$ on the interval $(-\infty, 1)$ \cite{Wakelin}. (This is the flow polynomial version of the analogous result for the chromatic polynomial \cite{T74}.) It follows that the sign of ${F}^{}_{H_3}({\phi}^{-2})$, and also the sign of ${F}^{}_{H_2}({\phi}^{-2})$ if $H_2$ is bridgeless, are opposite of that of ${F}^{}_{H_1}({\phi}^{-2})$. 
Rewrite  (\ref{ConjugateLinearRelation}) as 
$$
(-1)^{E-n+1}{F}^{}_{H_1}({\phi}^{-2})\,=  (-1)^{E-n} [ {\phi} \,  {F}^{}_{H_2}({\phi}^{-2})\, +\, {\phi}^2\, {F}^{}_{H_3}({\phi}^{-2}) ],
$$

where the first and the third terms are positive, and the second term is non-negative. Omitting the second term, one has
$$
(-1)^{E-n+1}{F}^{}_{H_1}({\phi}^{-2})\, \geq  (-1)^{E-n} [ {\phi}^2\, {F}^{}_{H_3}({\phi}^{-2}) ],
$$

To conclude the inductive step, note that that $H_3$ has two fewer vertices than $H_1$, and the inequality above gives a factor of ${\phi}^2$. The base of the induction is $F^{}_{K_4}=-{\phi}^2$.
\qed

It is interesting to compare the inequalities (\ref{Tutte inequality}) and (\ref{lower bound}). 
Conceptually the difference between the upper bound (\ref{Cauchy Schwartz}) at Beraha numbers and the lower 
bound (\ref{lower bound}) at ${\phi}^{-2}=(3-\sqrt{5})/2$ may be understood using the structure of the Temperley-Lieb algebra $TL^d$. As discussed in section \ref{background section}, the value $(3+\sqrt{5})/2$ of the chromatic polynomial corresponds to the parameter $d=2\, cos({\pi}/5)$ where $TL^d$ is positive semi-definite and the Cauchy-Schwarz inequality applies. (This is true for all Beraha numbers, see the proof of Theorem \ref{Cauchy Schwartz theorem}).
On the other hand the value $(3-\sqrt{5})/2$ corresponds to $d=2\, cos(3{\pi}/5)$.
The $4$th Jones-Wenzl projector is in the trace radical of $TL^d$ for both $d=2\, cos({\pi}/5)$ and $d=2\, cos(3{\pi}/5)$,  but the bilinear pairing on $TL^d$ at the latter value of $d$ has mixed signs.

\section{The golden identity and inequality for non-planar graphs}\label{inequality section}

The golden identity (\ref{golden identity}) for the chromatic 
polynomial of planar triangulations has an immediate counterpart stated dually for the flow polynomial of planar cubic graphs. Moreover, there are analogous identities
for the chromatic and flow polynomials relating their values at the Galois conjugate $(3-\sqrt{5})/2={\phi}^{-2}$ and $(5-\sqrt{5})/2$.
In this section we consider the validity of these relations for non-planar graphs. 

\subsection{A conjecture for the flow polynomial}

The Galois conjugate of the golden identity, stated for the flow polynomial of planar cubic  graphs, reads
$$  F^{}_G((5-\sqrt{5})/2) = (-{\phi}^{-1})^E\, F^{}_G((3-\sqrt{5})/{2})^2,$$
where $E$ is the number of edges of $G$.

Based on numerical evidence, we conjecture that in fact this identity characterizes planar cubic graphs:

\begin{conjecturesection} \label{golden conjecture} \sl
For any cubic bridgeless graph $G$,
\begin{equation} \label{conjugate inequality}
F^{}_G((5-\sqrt{5})/2) \geq (-{\phi}^{-1})^E\, F^{}_G((3-\sqrt{5})/{2})^2,
  \end{equation}
Moreover, $G$ is planar if and only if (\ref{conjugate inequality}) is an equality.
\end{conjecturesection}

Similarly, numerical evidence points to the inequality $ F^{}_G({\phi}+2) \leq {\phi}^E\, F^{}_G({\phi}+1)^2$ for any cubic graph $G$, with an equality if and only if $G$ is planar.

\subsection{The chromatic polynomial} 
\begin{figure}[ht]
\includegraphics[height=2.9cm]{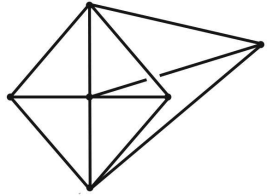}
\caption{A non-planar graph $G$ satisfying the golden identity (\ref{golden identity})}
\label{fig:NonPlanar}
\end{figure}

In contrast to Conjecture \ref{golden conjecture}, there are non-planar graphs satisfying the chromatic version (\ref{golden identity}) of the golden identity.
This section describes a specific example (see figure \ref{fig:NonPlanar}), and it suggests an explanation of the apparent difference with the behavior of the flow polynomial. 

The main point is that the {\em chromatic} Tutte linear relation (\ref{linear Tutte}) at ${\phi}+1$, 
 unlike its flow counterpart (\ref{3 term}), still holds when coupled with certain non-planar graphs.
(Note that in (\ref{linear Tutte})  there is an arbitrary number of edges  intersecting the boundary of the disk, so in fact this is an infinite collection of local relations for the chromatic polynomial of planar graphs.)

\begin{figure}[ht]
\includegraphics[height=2.9cm]{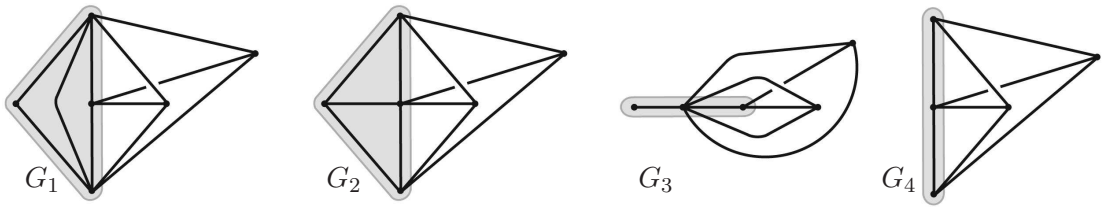}
{\small    \put(-414,10){$G_1$}
    \put(-300,10){$G_2$}
    \put(-182,10){$G_3$}
    \put(-90,10){$G_4$}  }
\caption{A non-planar graph $G_2$ satisfying the golden identity (\ref{golden identity}). Also 
shown are the graphs locally  related to $G_2$ as in figure \ref{fig:linear Tutte}.}
\label{fig:NonPlanarGolden}
\end{figure}

Consider the graph in figure  \ref{fig:NonPlanar}, obtained from $K_{3,3}$ by adding two edges. Endow it with the ribbon structure induced 
from the shown planar projection. A copy of this graph, denoted $G_2$, is also included in figure \ref{fig:NonPlanarGolden}. Consider the intersection of $G_2$ with the shaded region indicated in the figure as a particular example of the second graph in figure \ref{fig:linear Tutte}. The graphs $G_1, G_3, G_4$, shown below, are locally related to $G_2$ as in figure \ref{fig:linear Tutte}.
Note that $G_1, G_3, G_4$ are all planar, so they satisfy (\ref{golden identity}).

\begin{figure}[ht]
\includegraphics[height=2.4cm]{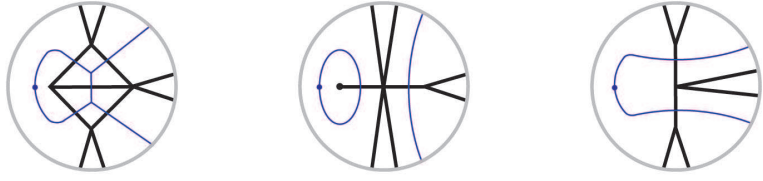}
\caption{Parts of the graphs $G_2,G_3,G_4$ (and their duals) in the shaded disks  in figure \ref{fig:NonPlanarGolden}.}
\label{fig:NonPlanar2}
\end{figure}

A key ingredient in the proof of the golden identity (\ref{golden identity}) in \cite{T2} (and also in \cite{FK})  is the linear identity (\ref{linear Tutte}) which holds 
as a consequence of the vanishing of the projector $P^{(4)}$ in the Temperley-Lieb algebra $TL^{\phi}$ (see section \ref{algebra map}).
This identity is intrinsic to planar graphs, and it certainly does not hold for general ribbon graphs. However we will now show that it still holds for {\em some} non-planar graphs. 

\begin{figure}[ht]
\includegraphics[height=2cm]{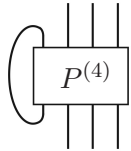}
{\small
    \put(-30,24){$P^{(4)}$}
}
\caption{The partial trace of the projector $P^{(4)}$ equals zero at $d={\phi}$ (and also at $d=-{\phi}^{-1}$).}
\label{fig:projector}
\end{figure}

In the relation (\ref{linear Tutte}) and figure \ref{fig:linear Tutte} there is an arbitrary number of edges leading to the boundary of the disk. In the special case in figure \ref{fig:NonPlanarGolden} the internal vertex on the left has 
{\em no} edges leading to the boundary. The map to the Temperley-Lieb algebra, described in section \ref{algebra map}, is defined in terms of dual trivalent graphs. The relevant dual graphs are shown in figure \ref{fig:NonPlanar2}. One checks that the result in the Temperley-Lieb algebra is the projector 
$P^{(4)}$ with two strands on the left connected, see figure \ref{fig:projector}. This is the partial trace of the projector $P^{(4)}$, and it equals zero 
in the Temperley-Lieb algebra at $d={\phi}$ \cite[Section 5.4]{FK}. Consequently, the identity (\ref{3 term}) and its chromatic dual:
\begin{equation} \label{projector equation}
{\phi}\, {\chi}^{}_{G_2}({\phi}+1)={\chi}^{}_{G_3}({\phi}+1)+(1-{\phi})\, {\chi}^{}_{G_4}({\phi}+1)
\end{equation}
hold for the graphs in figure \ref{fig:NonPlanarGolden}. (More generally, the identity still holds for an arbitrary number of edges leading to the 
boundary of the disk at the other three vertices in figure \ref{fig:NonPlanar2}, and these local pictures may be paired with an arbitrary non-planar graph.)

The chromatic polynomial of these graphs  satisfies the contraction-deletion rule pictured in figure~\ref{chromatic relations}:
\begin{equation} \label{eq1}
{\chi}^{}_{G_1}({\phi}+2)+{\chi}^{}_{G_4}({\phi}+2)={\chi}^{}_{G_2}({\phi}+2)+{\chi}^{}_{G_3}({\phi}+2).
\end{equation} 
As shown in the proof
of lemma 3.1 in \cite{FK}, (\ref{projector equation}) yields
\begin{equation} \label{eq2}
{\phi}^3\, ({\chi}^{}_{G_1}({\phi}+1))^2+({\chi}^{}_{G_4}({\phi}+1))^2={\phi}^3\, ({\chi}^{}_{G_2}({\phi}+1))^2+({\chi}^{}_{G_3}({\phi}+1))^2.
\end{equation}
Equations (\ref{eq1}), (\ref{eq2}), combined with the fact that three of the graphs satisfy the golden identity, imply that
the remaining graph, $G_2$, also obeys (\ref{golden identity}).

The argument above suggests a strategy for generating families of non-planar graphs satisfying the chromatic golden identity.
It would be interesting to find a characterization of all such graphs, and in particular to see
whether the argument using the 4th Jones-Wenzl projector considered here generates all of them.
Another open question is whether there exist $4$-connected non-planar graphs with the number of edges and vertices satisfying  $E=3V-6$ (analogous to the case of a planar triangulation) and satisfying the golden identity  (\ref{golden identity}).

{\bf Acknowledgments}.  
The authors would like to thank the IHES for hospitality and support (NSF grant 1002477). 
Ian Agol was partially supported by NSF grant DMS-1406301, and by a Simons Investigator grant.
Vyacheslav Krushkal was partially supported by NSF grant DMS-1309178.

\end{document}